\documentclass{article}
\usepackage{epsfig}
\usepackage{latexsym}
\usepackage{amsfonts}
\usepackage{amssymb}
\def\area{{\mathrm{area}\,}}
\def\C{{\mathbf{C}}}
\def\bC{{\mathbf{\overline{C}}}}
\title{Transcendental meromorphic functions
with three singular values}
\author{A. Eremenko\thanks{Supported by NSF grants DMS-0100512,
DMS-0244421,
and by the
Humboldt Foundation.}}
\date{February 27, 2004}
\begin{document}
\maketitle
\begin{abstract}
Every transcendental meromorphic function $f$
in the plane, which has only three critical values,
satisfies
$$\liminf_{r\to\infty}\frac{T(r,f)}{\log^2r}\geq
\frac{\sqrt{3}}{2\pi},$$
and this estimate is best possible.

\end{abstract}

A {\em singular value} of a meromorphic function $f$ in the
plane $\C$ is, by definition,
a critical value or an asymptotic value.
If we denote the closure
of the set of singular values by $S$, then
$$f:\C\backslash f^{-1}(S)\to\bC\backslash S$$
is a covering map. Meromorphic functions with
finitely many
singular values play an important role in
value distribution
theory (see, for example, \cite{GO,Teich,Wittich}, and
Jim Langley's papers on the distribution of
values of derivatives),
as well as in holomorphic dynamics
\cite{Berg,EL}.

In this paper, ``meromorphic function''
will always mean a {\em transcendental\footnote{\rm On algebraic
functions with three critical values see
\cite{G1,G2}.}
meromorphic function in the plane}, unless some other region is specified. 

Langley \cite{Lang1,Lang3} discovered that there exists a relation between the number of
singular values of a meromorphic function $f$
and the growth of the Nevanlinna characteristic $T(r,f)$.
In \cite{Lang3} he proved that all meromorphic functions $f$ with
finitely many singular values satisfy
$$\liminf_{r\to\infty}\frac{T(r,f)}{\log^2r}>0.$$
On the other hand, he showed
in \cite{Lang2}
that for every $\epsilon>0$
there exists a meromorphic function
$f$ with four singular values such that
$$\limsup_{r\to\infty}\frac{T(r,f)}{\log^2r}<\epsilon.$$
Concerning meromorphic functions with three singular values,
Langley proved in \cite{Lang2} that they satisfy
\begin{equation}
\label{1}
\liminf_{r\to\infty}\frac{T(r,f)}{\log^2r}\geq c,
\end{equation}
where $c$ is an absolute constant.

In this paper, the precise value of this
constant is found. I thank Walter Bergweiler
who brought \cite{Lang2} to my attention and
suggested the extremal problem which the following
theorem solves.
\vspace{.1in}

\noindent
{\bf Theorem} {\em Let $f$ be a meromorphic function 
with at most three singular values. 
Then $(\ref{1})$ holds with
$c=\sqrt{3}/(2\pi)$, and there exists a meromorphic function $f_0$ 
with three singular values, such that}
\begin{equation}\label{Tre}
T(r,f_0)/\log^2r\to\sqrt{3}/(2\pi)\quad\mbox{as}\quad r\to\infty.
\end{equation}

\noindent
{\bf Remarks}
\vspace{.1in}

1. If $f$ has finitely many singular values and at least one of them
is an asymptotic value, then
\begin{equation}
\label{2}
\liminf_{r\to\infty}\frac{T(r,f)}{\sqrt{r}}>0;
\end{equation}
in particular this holds for all entire
functions with finitely
many singular values.
We sketch a proof of this, different from
that mentioned in \cite{Lang2}.
If $a$ is an asymptotic value, and there are
no other singular values in an $\epsilon$-neighborhood
of $a$, then 
one of the components
of the set $\{ z:|f(z)-a|<\epsilon\}$ is
an unbounded region $D$ whose boundary consists
of one simple curve, and $f(z)\neq a$ in this region.
Applying the standard growth estimates 
to the harmonic function
$\log|f(z)-a|^{-1}$ in $D$ we conclude that (\ref{2}) holds.
\vspace{.1in}

2. {\em Let $f$ be a meromorphic function with at most three critical values.
Then the conclusion of the Theorem holds.}
\vspace{.1in}

Indeed, we can assume that $f$ has finite
lower order (otherwise there is nothing to prove).
For such functions with finitely many critical values, the
set of asymptotic values is also finite. 
This was proved in \cite{BE} for
functions of finite order and extended in
\cite{H} to functions of
finite lower order.
If there are asymptotic values, we apply
Remark 1, if there are none, we apply the Theorem.

A similar improvement can be made in the results of Langley mentioned above.
\vspace{.1in}
 
3. All meromorphic functions with two singular values
are of the form $L\circ\exp$,
where $L$ is a fractional-linear
transformation.
\vspace{.1in}

4. Let $\mu$ be a probability measure on the Riemann sphere $\bC$,
and $\nu=f^*\mu$ the pull-back of $\mu$. We denote
$$A_\mu(r)=\nu\left(\{ z:|z|\leq r\}\right),\quad r>0.$$
If $\mu$ is the normalized spherical area, then $A_\mu(r)\equiv A(r)$
is the average number of sheets of 
the map $f:\{ z:|z|\leq r\}\to\bC$.
The Nevanlinna characteristic satisfies
$$T(r,f)=
\int_e^r\frac{A(t)}{t}dt+O(\log r),\quad r\to\infty.$$
For an arbitrary probability measure $\mu$ on the sphere, we have
$$\int_e^r\frac{A_\mu(t)}{t}dt\leq T(r,f)+O(1),$$
which is a consequence of the First Fundamental Theorem of Nevanlinna,
\cite[VI,\S 4]{Nev}.
Thus, to prove our theorem, it is sufficient to show
that
\begin{equation}
\label{3}
A_\mu(t)\geq\frac{\sqrt{3}}{\pi}\log t+O(1),\quad t\to\infty,
\end{equation}
for some probability measure $\mu$.
\vspace{.1in}

{\em Proof of the Theorem}.
\vspace{.1in}

Let $F$ be a connected oriented surface of finite
topological type, possibly with boundary.
A {\em triangular net} on $F$
is a locally finite covering
of $F$ by closed sets $T$ called {\em triangles}, such that
\vspace{.1in}

(i) Each triangle $T$ is a closed Jordan region (homeomorphic
to a closed unit disc) with three marked
distinct boundary points called {\em vertices}.
A closed boundary arc between two adjacent vertices
is called an {\em edge} of $T$. 

(ii) The intersection of two triangles is either empty, or 
a union of common edges and common vertices.

(iii) Triangles are divided into two classes, white and
black, so that any two triangles with a common
edge are of different colors.

(iv) Vertices are labeled by the letters $A,B$ and $C$,
so that for each triangle, all three labels are
present on its boundary.
Furthermore, the cyclic order on the oriented boundary of
a white triangle is $(A,B,C)$, and opposite on the oriented boundary
of a black triangle.
\vspace{.1in}

Suppose that a covering of $F$ by triangles satisfying (i) and (ii)
is given, and there exists a labeling of vertices satisfying (iv).
Then such a labeling is uniquely defined by a choice of labels
on the three vertices of one triangle. The colors of triangles
are evidently determined by the labeling of vertices.
\vspace{.1in}

Let $f$ be a meromorphic function satisfying the assumptions
of the Theorem. In view of Remarks 1--3,
we assume without loss of generality that $f$ has no asymptotic values
and exactly three critical values, $A,B$ and $C$. Composing $f$ with
a fractional linear transformation we may assume 
that $A,B$ and $C$ are real, and $A<B<C$.

Consider the triangular net $N_0$ on the Riemann sphere,
which consists of
two triangles, the closed upper half-plane (white) and the
closed lower half-plane (black), and the vertices
$A,B$ and $C$. 

We construct the $f$-preimage of this net, which will be called 
$N$. The vertices of $N$ are preimages of the vertices of $N_0$,
and they are labeled according to their images.
The triangles of $N$ are the closures of the components of
the preimages of the open
upper and lower half-planes, and the edges of $N$ are
defined in the evident way. Each triangle in $N$ is assigned the
same color as its image in $N_0$.

Our assumptions about $f$ imply that $N$ is a triangular net
in the plane $\C$.

Choose one triangle $T_0\in N$ 
and remove its interior
${\mathrm{int}\,}T_0$ from the plane.
We obtain a bordered surface 
$D=\C\backslash {\mathrm{int}\,}T_0$ which is
homeomorphic to a closed
semi-infinite cylinder.

Let $G\subset D$ be a compact closed region
homeomorphic to a closed ring,
which separates $\partial T_0$ from $\infty$.
We are going to estimate
from below
the {\em extremal length}\footnote{
See, for example, \cite{Ahlfors} for a definition
and simplest properties.}
$\lambda$
of the family of all closed non-contractible curves 
in $G$. For this purpose we construct
a conformal metric $\rho$ in
$D$. In this metric, each triangle
of the net $N\backslash T_0$ will be isometric to a
Euclidean equilateral triangle $\Delta$ with sidelength $1$.

First we define a flat conformal metric $\rho_0$ on
$\bC\backslash\{ A,B,C\}$.
Let $g$ be the conformal map of $\Delta$ onto the upper half-plane
sending the vertices of $\Delta$ to $A,B,C$.
(An explicit expression of $g$ will be given at the end of the
paper). The metric $\rho_0$ is defined in
the upper half-plane by the length element
$ds=|(g^{-1})'(z)|\,|dz|$, so that $g$ becomes an isometry
from $\Delta$ with the Euclidean metric to the upper
half-plane with the metric $\rho_0$.
Using the Symmetry Principle, we extend $\rho_0$ to
$\bC\backslash\{ A,B,C\}$. 

The Riemann sphere with the metric $\rho_0$ can
be visualized as a ``two-sided triangle $\Delta$''.
The area of the sphere with respect to $\rho_0$ is
$\sqrt{3}/2$.

Now we define $\rho=f^*\rho_0$, the pull-back of $\rho_0$
via $f$.

The metrics $\rho$ and $\rho_0$ 
have isolated singularities at the vertices,
but this does not cause any problems.

Now we estimate the $\rho$-length of 
non-contractible
curves in $D$ from below. 
\vspace{.1in}

\noindent
{\bf Lemma} {\em Let $D'$ be a surface homeomorphic to $\{1\leq|z|\leq 2\}$,
equipped
with a triangular
net and 
an intrinsic metric\footnote{A metric is called intrinsic if
the distance between any two points equals the infimum of lengths of
curves
connecting these points.} such that every triangle of the net is isometric 
to $\Delta$. Then the length of every closed non-contractible curve
in $D'$ is at least $\sqrt{3}$.}
\begin{center}
\begin{picture}(0,0)%
\includegraphics{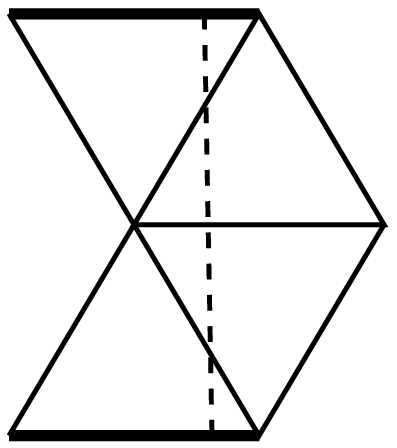}%
\end{picture}%
\setlength{\unitlength}{1973sp}%
\begingroup\makeatletter\ifx\SetFigFont\undefined%
\gdef\SetFigFont#1#2#3#4#5{%
  \reset@font\fontsize{#1}{#2pt}%
  \fontfamily{#3}\fontseries{#4}\fontshape{#5}%
  \selectfont}%
\fi\endgroup%
\begin{picture}(3730,4255)(3515,-6105)
\put(4576,-2836){\makebox(0,0)[lb]{\smash{\SetFigFont{12}{14.4}{\rmdefault}{\mddefault}{\updefault}{$T_3$}%
}}}
\put(5701,-3586){\makebox(0,0)[lb]{\smash{\SetFigFont{12}{14.4}{\rmdefault}{\mddefault}{\updefault}{$T_1$}%
}}}
\put(4576,-5386){\makebox(0,0)[lb]{\smash{\SetFigFont{12}{14.4}{\rmdefault}{\mddefault}{\updefault}{$T_4$}%
}}}
\put(5701,-4936){\makebox(0,0)[lb]{\smash{\SetFigFont{12}{14.4}{\rmdefault}{\mddefault}{\updefault}{$T_2$}%
}}}
\end{picture}

\vspace{.1in}

Figure 1. An extremal configuration.
 Bold segments are identified. An extremal curve is the broken line.
\end{center}
\vspace{.1in}

The idea of the following proof, which is simpler than
the original one, was suggested by Mario Bonk. 
\vspace{.1in}

{\em Proof}. Let $\gamma$ be a shortest non-contractible
curve in $D'$; evidently such a curve exists.
Then $\gamma$ is homeomorphic to a circle,
because otherwise we could remove extra loops and
shorten~$\gamma$.

Suppose first that $\gamma$
passes through a vertex $v$ of the net, that is $\gamma(t_0)=v$.
Let $F$ be the interior of the union of all
triangles of the net that have $v$ as a common vertex. Then $F$
is a simply connected region, and to be 
non-contractible,
our curve has to pass through a point $w\in\partial_D F=
\partial F\cap{\mathrm{int}\,} D$.
Then two arcs of $\gamma$ from $v$ to $w$ have
lengths at least $\sqrt{3}/2$ each, which proves the Lemma in
this case.
\vspace{.1in}

\noindent
{\em From now on we assume that $\gamma$ does not pass through
the vertices.}
\vspace{.1in}

\noindent
As our metric $\rho$ is flat away from
the vertices, every point $w\in D'\backslash \{\mbox{vertices}\}$
has a neighborhood $W$
which is isometric to a region $V$ in the plane with the standard
(intrinsic) metric. This isometry $\phi:W\to V$ is called
the ``developing map''.
It has an analytic continuation to every simply
connected region in $D'\backslash\{{\mathrm{vertices}}\}$
that contains $w$, and also an analytic continuation
along any curve in $D'$ which does not
pass through the vertices. The image of any arc of our shortest curve $\gamma$
under the developing map is an interval of a straight line in the plane.

Our curve $\gamma$ evidently intersects some edge.
Let $v$ be a point of intersection with an edge $e$,
and choose a parametrization 
$\gamma:[0,1]\to D'$ such that $\gamma(0)=\gamma(1)=v$. Let $\phi$
be a germ of the developing map at $v$, and $\phi_1$ the result
of its analytic continuation along $\gamma$. Let $T_1,\ldots,T_n$
be the sequence of triangles in $D'$ visited by $\gamma$,
enumerated in the natural
order. Then the branches of the developing map are defined in $T_j$
by analytic continuation of the germ $\phi$ along $\gamma$. Let
$\Delta_1,\ldots,\Delta_n$ be the images of the triangles
$T_1,\ldots,T_n$ under these branches.

Reflections in the sides of $\Delta_1$ generate a group $\Gamma$ of
isometries of the plane. It is clear
that $\Delta_1,\ldots\Delta_n$ can be obtained
by applying some elements of this group to
$\Delta_1$. The full orbit of $\Delta_1$ under $\Gamma$ forms the hexagonal
tiling of the plane by equilateral triangles of sidelength $1$.
We label the vertices
of $\Delta_1$ similarly to the corresponding labels of $T_1$.
This uniquely defines the labeling of vertices of the hexagonal
tiling of the plane, thus making it into a triangular net in $\C$.
The branches of the developing map 
$T_j\mapsto\Delta_j$ preserve the labels of vertices.  

The edges of our nets have orientations induced by the
cyclic order on the labels of vertices.
Let us consider the edges $e^*=\phi(e)$ and $e^*_1=\phi_1(e)$.
We claim that the oriented edge $e^*_1$ can be obtained from
the oriented edge $e^*$ by a translation from $\Gamma$ which is not
the identity map.

Let $\alpha$ be the angle, counted anti-clockwise
from the positive direction of $e$ to the tangent vector 
$\gamma'(0)$ at $v$. The image of $\gamma$ under the analytic
continuation of the developing map is a (non-degenerate)
straight line segment which makes
the
same angle $\alpha$ with $e^*$ and $e_1^*$. We conclude
that $e^*$ and $e_1^*$ have the same
direction but do not coincide.
As $e^*$ and $e^*_1$ are edges of
the hexagonal triangular net, we easily conclude that $e^*_1$
can be obtained from $e^*$ by a translation in $\Gamma$.
This proves our claim.

Now, he shortest translation in $\Gamma$ has
magnitude $\sqrt{3}$, and this completes the proof of the Lemma.
\hfill$\Box$
\vspace{.1in}

As a corollary we obtain that the extremal length $\lambda$
of the set of non-contractible curves in any compact ring
$G\subset D$
separating $\partial T_0$ from $\infty$ in $D$ satisfies
\begin{equation}
\label{length}
\lambda\geq\frac{3}{\area(G)},
\end{equation}
where the area corresponds to the metric $\rho$.

To complete the proof, we consider the set
$$G(t)=\{ z:|z|\leq t\}\backslash{\mathrm{int}\,} T_0.$$
This set is homeomorphic to a ring if $t$ is large enough.
The extremal length $1/\lambda$ of the family of curves
connecting the boundary components of the ring $G(t)$
is $(2\pi)^{-1}\log t+O(1)$ as $t\to\infty$.
According to (\ref{length}), the area of this
ring with respect to the metric $\rho$
is at least $$\frac{3}{\lambda}=\frac{3}{2\pi}\log t+O(1),\quad t\to\infty.$$
The $\rho_0$-area of the Riemann
sphere is $\sqrt{3}/2$. Taking $\mu$ to be the $\rho_0$-area divided by
$\sqrt{3}/2$, we obtain  $A_\mu(t)\geq
(\sqrt{3}/\pi)\log t+O(1)$, which is (\ref{3}).
\hfill$\Box$
\vspace{.1in}

\noindent
{\bf Example}
\vspace{.1in}

Consider the equilateral triangle
$\Delta\subset\C$ with vertices $0,i$ and
$(\sqrt{3}+i)/2$. Let $g$ be a conformal
map of this triangle onto the right half-plane, 
sending the vertices to $\infty,ia,-ia$, where $a>0$. By reflection,
$g$ extends to a meromorphic function in $\C$
with no asymptotic values and three critical values,
$ia,-ia$ and $\infty$. All preimages of the
critical values are critical points of order $3$.

This function $g$ is doubly periodic and its shortest
period is $\sqrt{3}$. It can be expressed in terms
of the Weierstrass function of an equiharmonic lattice,
see \cite{AS,Ford}.
If we choose
$$a=k^3,\quad\mbox{where}\quad k=\frac{\Gamma^3(1/3)}{2\pi\sqrt{3}},$$
then $g=\wp'$ where $\wp$ is the Weierstrass function
with periods $\sqrt{3}$ and $\sqrt{3}e^{2\pi i/3}$.
The Riemannian metric $\rho_0$ in the proof of the Theorem
corresponds to the length element $|(g^{-1})^\prime(w)|\,|dw|.$

The function
$$f_1(z)=
g\left(\frac{\sqrt{3}}{2\pi i}\log z\right),$$
is evidently meromorphic in $\C^*=\C\backslash\{0\}$
and has three critical values.
Simple calculation shows that it satisfies
$$A(r,f_1)=\frac{\sqrt{3}}{\pi}\log r+O(1),\quad r\to\infty.$$
Now we modify $f_1$ to obtain a function meromorphic
in $\C$. Consider the integral
$$I(z)=\int_0^z\frac{d\zeta}{\zeta^{2/3}(1-\zeta)^{1/3}},$$
Using the positive value of the cubic root,
and integrating over $[0,1]$ we obtain
$$I(1)=2\pi/\sqrt{3}.$$
Using this, one can easily verify that 
$$f_0(z)=g\left(\frac{\sqrt{3}}{2\pi i}I(z)\right)$$
is meromorphic in the plane. This function $f_0$
has three critical values,
no asymptotic values, and satisfies (\ref{Tre}), as a branch of $I(z)$
near infinity
has the same asymptotic behavior as the logarithm.
\vspace{.1in}

We also sketch a purely geometric construction of $f$,
based on the Uniformization theorem.

Consider the closed region $K_0$ in $\C$
made of four equilateral triangles
as in Figure 1, and put
$$K=\cup_{j=1}^\infty (K_0+j),$$
the union of translates of $K_0$ by positive integers.
Then $K$ is a closed half-strip in the plane, and we define a Riemann surface
$P$ by identifying pairs of points with equal
$x$-coordinates on the two horizontal boundary rays of $K$. 
The surface $P$ is homeomorphic to a semi-infinite closed
cylinder, and its boundary consists of two edges.
We ``patch the hole'' by adding a $2$-gon consisting of two
triangles with two common edges.

Thus we obtain a surface $P'$ homeomorphic to the plane and covered
by triangles.
It is easy to see that one can label the vertices of triangles
in $P'$ to obtain a triangular net. So we obtain a triangular net $N$
in the plane.

Every triangular net in the plane defines a ramified covering $h:\C\to\bC$
in the following way. Choose a triangular net $N_0$ of $\bC$ with
two triangles and three vertices. First define $h$ on the vertices,
by sending each vertex of $N$ to the similarly labeled vertex of $N_0$. Then
extend $h$ to the edges, so that the restriction to
each edge is a homeomorphism, and the extended map
is continuous on the $1$-skeleton of $N$.
Finally, extend $h$ to the interiors of triangles, so that the extended
map is a homeomorphism on each triangle. 

By the Uniformization theorem, there exists a homeomorphism $\psi:U\to\C$,
where $U$ is either a disc or the plane, such that
$h\circ\phi$ is a meromorphic function.

To show that $U=\C$, and to prove (\ref{Tre}) one can make extremal length
estimates
using the Euclidean metric on $P'$ such that all
triangles of the net are equilateral with side length $1$.
We omit the details which are similar to those in the proof of the Theorem.

The author thanks A.A. Go'ldberg and V. S. Azarin for a useful discussion
of this paper.

\vspace{.1in}

Purdue University

West Lafayette, IN 47907

{\em eremenko@math.purdue.edu}

\begin{thebibliography}{1}
\bibitem{AS} M. Abramowitz and I. Stegun, {\em Handbook of mathematical
functions
with formulas, graphs, and mathematical tables,} Natl. Bureau of Standards,
Washington DC, 1964.
\bibitem{Ahlfors} L. Ahlfors, {\em Conformal invariants,}
McGraw-Hill, NY, 1973.
\bibitem{G1} G. Belyi,
Galois extensions of a maximal cyclotomic field, 
{\em Izv. Akad. Nauk SSSR Ser. Mat.} 43 (1979) 267--276, 479. 
\bibitem{Berg} W. Bergweiler, Iteration of meromorphic
functions, {\em Bull. Amer. Math. Saoc.}, 29 (1993) 151--188.
\bibitem{BE} W. Bergweiler and A. Eremenko, On the singularities
of the inverse to a meromorphic function of finite order,
{\em Rev. Mat. Iberoamericana,} 11 (1995) 355--373.
\bibitem{EL} A. Eremenko and M. Lyubich, Dynamics of
some classes of entire functions,
{\em Ann. Inst. Fourier,} 42 (1992) 1--32.
\bibitem{Ford} L. Ford, {\em Automorphic functions,} Chelsea, NY, 1951.
\bibitem{GO} A.A. Goldberg and I.V. Ostrovskii, {\em Distribution of
values of meromorphic functions,} Moscow, Nauka, 1970 (Russian)
\bibitem{G2} A. Grothendieck, 
Esquise d'un programme,  in: L. Schneps, ed., 
{\em Geometric Galois actions,} 1, 5--48, 
Cambridge Univ. Press, Cambridge, 1997. 
\bibitem{H}
J.D. Hinchliffe, 
The Bergweiler-Eremenko theorem
for finite lower order, 
{\em Results Math.,} 43 (2003), no. 1-2, 121--128.
\bibitem{Lang1} 
J.K. Langley, On the multiple points of certain meromorphic functions,
{\em Proc. Amer. Math. Soc.,} 123, no. 6 (1995), 1787--1795.
\bibitem{Lang3} J.K. Langley, On differential polynomials,
fixpoints and critical values of meromorphic functions,
{\em Result. Math.,} 35 (1999) 284--309.
\bibitem{Lang2} J.K. Langley, Critical values of slowly
growing meromorphic functions, to appear in Computational Methods and
Function Theory.
\bibitem{Nev} R. Nevanlinna, {\em Eindeutige analytische Funktionen,}
Springer, Berlin, 1953.
\bibitem{Teich} O. Teichm\"uller, Eine Umkehrung des zweiten Hauptsatzes
der Wertverteilungslehre, {\em Deutsche Math.,} 2 (1937) 96--107;
{\em Ges. Abh.,} Springer, Berlin, 1982, 158--169.
\bibitem{Wittich} H. Wittich, {\em Neuere Untersuchungen
\"uber eindeutige analytische Funktionen,}
Springer, Berlin, 1955.
\end{thebibliography}
\end{document}